\documentclass{amsart}
\usepackage{amssymb}
\usepackage{amsthm}

\theoremstyle{plain}
\newtheorem{thm}{Theorem}[section]
\newtheorem{lem}[thm]{Lemma}
\newtheorem{cor}[thm]{Corollary}
\newtheorem{prop}[thm]{Proposition}

\theoremstyle{definition}
\newtheorem{dfn}[thm]{Definition}
\newtheorem{Q}[thm]{Question}

\newcommand{\Fp}{\mathbb{F}_p}

\begin{document}

\title{The isomorphism problem for classes of computable fields}
\author{Wesley Calvert} 
\address{Department of Mathematics \\ 255 Hurley Hall \\ University of Notre Dame \\ Notre Dame, Indiana 46556}
\email{wcalvert@nd.edu}
\thanks{The author was partially supported by NSF Grants DMS 9970452 and DMS 0139626. The author wishes to thank J.\ 
F.\ Knight for many helpful comments on this article.} 
\begin{abstract} Theories of classification distinguish classes with some good structure theorem from those for which none is possible. Some classes (dense linear orders,
for instance) are non-classifiable in general, but are classifiable when we consider only countable members. This paper explores such a notion for classes of computable
structures by working out several examples. One motivation is to see whether some classes whose set of countable members is very complex become classifiable when we
consider only computable members.

We follow recent work by Goncharov and Knight in using the degree of the isomorphism problem for a class to
distinguish classifiable classes from non-classifiable. For arbitrary fields --- even real closed fields --- we show
that the isomorphism problem is $\Sigma^1_1$ complete (the maximum possible), and for others we show that it is of
relatively low complexity.  We show that the isomorphism problem for algebraically closed fields, Archimedean real
closed fields, or vector spaces is $\Pi^0_3$ complete.
\end{abstract}
\maketitle

\section{Introduction}

We will consider a notion of ``classification" for classes of computable structures. For some classes, there is a
``classification," or ``structure theorem" of some kind. For instance, the classification of algebraically closed
fields states that a single cardinal (the transcendence degree) completely determines the structure up to isomorphism.
For other classes (graphs, for example, or arbitrary groups) such a result would be surprising, and when we introduce
the necessary rigor we can prove that there is none to be found. They simply have more diversity than any structure
theorem could describe.

A theory of classification should tell us which classes fall into which of these two categories. Such a theory, 
originating in the work of Shelah, has
long been in use for elementary classes (see \cite{shelah,hodges}). There is
also a well-developed theory for classes of countable structures, which uses Borel reducibility.  This theory was 
developed by Friedman and Stanley \cite{fs}, Becker, Hjorth, Kechris, Louveau, and others \cite{bkbook,hjkelo}. 

Shelah's theory of classifications considers structures of arbitrary cardinality. The Borel reducibility notion
focuses on countable structures only. This difference is exemplified by the theory of dense linear orders, which is
$\aleph_0$-categorical, but whose class of models is non-classifiable in the Shelah sense because for any uncountable
$\kappa$ there are $2^\kappa$ non-isomorphic models of size $\kappa$.  However, in some classes which are
non-classifiable in the Shelah sense, the complexity is in some way so intrinsic to the theory that it shows up even
in the countable case, in the form of high Borel complexity.

We work only with structures which have for their universe a subset of $\omega$, and identify structures with their
atomic diagrams. Thus, for instance, a structure is computable if and only if its atomic diagram is computable as a
set of G\"{o}del numbers of sentences. Alternatively, we could use the quantifier-free diagram instead of the atomic
diagram. Similarly, a structure is associated with the index of a Turing machine which enumerates its atomic diagram
(assuming its universe is computable). In this paper I will write $\mathcal{A}_a$ for the computable structure with
atomic diagram $W_a$ and will consider only classes $K$ of structures which have only computable members.

\begin{Q} Are there some classes where this high complexity is apparent even in the class of computable
models?\end{Q}

The answer is not immediately apparent. Any class of countable structures with high Borel complexity has
$2^{\aleph_0}$ members up to isomorphism. It might seem that if we restrict to only countably many -- and at that the
most tangible members, the computable ones -- we might have required so much regularity that the enormous complexity
we saw before would be forbidden. It often happens that there is some structural characterization of which members of
a class admit a computable structure, and perhaps such results (known and unknown) should mean that the set of
computable members of some class cannot be too complicated.

Surprising or not, it turns out that we can still distinguish at the computable level between some ``very complicated"
classes and some ``quite simple" ones. It is easy to tell whether two algebraically closed fields are isomorphic, and
it is difficult to do the same for arbitrary fields, even if we only consider the computable models of each theory.

We are now prepared to make precise what I have meant by the terms ``simple" and ``complicated." The following
formalization was recently proposed by Goncharov and Knight \cite{gk}.

\begin{dfn}The isomorphism problem, denoted $E(K)$, is the set \[\{ (a,b) | \mathcal{A}_a , \mathcal{A}_b \in
K\mbox{, and }\mathcal{A}_a \simeq \mathcal{A}_b\}\]\end{dfn}

If the set of indices for computable members of $K$, denoted $I(K)$, is hyperarithmetical, then $E(K)$ is 
$\Sigma^1_1$.
Intuitively, in the worst case, where $E(K)$ is $\Sigma^1_1$ complete, the easiest way to say that two members of $K$
are isomorphic is to say, ``There exists a function which is an isomorphism between them." Often there are easier ways
to check isomorphism, such as counting basis elements of vector spaces. Such a ``shortcut" is a classification. As is
discussed more fully in \cite{gk} and in section 3 of the present paper, this notion is closely related to more common
understandings of classification, such as the production of a ``nice" list of isomorphism types.

Several classes are well-known to have maximally complicated isomorphism problems.  The following theorem summarizes
several classical results.  Proofs may be found in articles by Rabin and Scott \cite{rabsco}, Goncharov and Knight
\cite{gk}, Morozov \cite{moro}, and Nies \cite{nies}.
\begin{thm} If $K$ is the set of computable members of any of the following classes, then $E(K)$ is $\Sigma^1_1$ 
complete:
\begin{enumerate}
\item Undirected graphs
\item Linear orders
\item Trees
\item Boolean algebras
\item Abelian $p$-groups
\end{enumerate}
\end{thm}

The following additions to the list follow easily from recent work by Hirschfeldt, Khoussainov, Shore, and Slinko
\cite{hkss}.
\begin{thm} [Hirschfeldt -- Khoussainov -- Shore -- Slinko] If $K$ is the set of computable members of any of the 
following classes, then $E(K)$ is $\Sigma^1_1$ complete:
\begin{enumerate}
\item [6.] Rings
\item [7.] Distributive lattices
\item [8.] Nilpotent groups
\item [9.] Semigroups
\end{enumerate}
\end{thm}

This part of the paper will focus on calculating the complexity of the isomorphism problem for classes of fields. We
will make use of the example of undirected graphs, so in section 2 this example is worked out.  Section 3 contains a
proof that for the class of arbitrary computable fields, the isomorphism problem is maximally complicated. In section
4, we concentrate on algebraically closed fields, where the isomorphism problem is quite simple, just as we would
expect. The example of vector spaces is a comfortable warm up for algebraically closed fields, so it is also included
in this section.  Section 5 treats real closed fields, both Archimedean (where the isomorphism problem is rather
simple) and arbitrary (where it is $\sum^1_1$ complete). In a later paper, I will consider certain classes of
Abelian groups \cite{ecl2}. 

The goal of this paper is twofold. On one hand, the isomorphism problem gives us considerable insight into these
classes and the great diversity possible among, for example, computable fields. On the other hand, these classes serve
as benchmarks to show that the complexity of the isomorphism problem gives the ``right" answers to distinguish between
classifiable and non-classifiable classes.

\section{Undirected Graphs}

An older result, of which proofs are given in \cite{moro} and \cite{gk},
shows that the isomorphism problem for the class of computable directed
graphs is $\sum^1_1$ complete. Friedman and Stanley \cite{fs} state as
well-known the fact that the isomorphism problem for countable undirected
graphs is Borel complete. It is also known that the isomorphism for
undirected graphs is $\sum^1_1$ complete. Morozov \cite{moro}, Nies
\cite{nies}, and Rabin and Scott \cite{rabsco} each give a proof. Below,
we give a slightly simplified version of Nies's proof. The idea of the
construction is that if there is a directed edge from $n$ to $m$ in our
directed graph, we should have a connection between them in the undirected
graph which is labeled in a particular way: $n$ is connected to a triangle,
which is connected to a pentagon, which is connected to $m$. Since we can
distinguish the triangle from the pentagon, we can determine the direction
of the connection.

The difficulty is in distinguishing which points are the actual ``vertices" and which are the edge labels. This will be done by having an ``Archimedean point" to which all of the ``vertices" are connected, but none of the edge-labeling points. If the original directed graph was completely disconnected, then the Archimedean point will be the only one connected to infinitely many others. Otherwise, it will be the only point $x$ where the pattern $x$-$y$-triangle-pentagon-$z$-$x$ occurs. The ``vertex" points are exactly those points which are directly connected to the Archimedean point. I do not mean to say that this information may be obtained effectively; only that it is isomorphism invariant.

\begin{thm} When $K$ is the class of computable undirected graphs, $E(K)$ is $\sum^1_1$ complete.\end{thm}

\begin{proof} The class of graphs (irreflexive, symmetric binary relations) is characterized by $\prod^0_1$ axioms, so $I(K)$ is $\prod^0_2$ (identifying a set as the atomic diagram of some structure is $\prod^0_2$ complete). Given a directed graph $T$, we will produce an undirected graph $G$, as described above. This construction will be uniformly computable in $T$, and will be injective on isomorphism types. We will use $1$ for the Archimedean point, and we start with $D(G)_{-1} = P_{-1} = C_{-1} = V_{-1} = \nu_{-1} = \emptyset$.

Step $6s$: We will go a bit farther on connecting the Archimedean point to the vertices. In particular, for the first $x \in V_{s-1}$, if any, such that $1Gx \notin D(G)_{6s-1}$ set $D(G)_{6s} = D(G)_{6s-1} \cup \{1 G x\}$. If there is no such $x$, set $D(G)_{6s} = D(G)_{6s-1}$.

Step $6s+1$: In this step, we will check a connection from $T$, and if 
necessary, represent it in $G$.  Let $g(s) = (n_s, m_s)$ be a pairing 
function. First we should check whether $n_s$ and
$m_s$ are already represented as vertices (i.e. if $n_s, m_s \in dom(\nu_{s-1})$). If not, suppose $n_s$ is not represented. Then find the least $y$ not occurring in the
diagram so far, and set $V_s = V_{s-1} \cup \{y\}$ and $\nu_s = \nu_{s-1} \cup \{(n_s, y)\}$. If necessary, act similarly for $m_s$. Now since $T$ is computable, we can
check whether $n_s T m_s$. If we find that there is such a connection, and if $(n_s, m_s) \notin C_{3s-1}$ we find the least $x_0, \dots, x_7$ which do not occur in the
(finite) diagram so far, and which are not equal to $1$. We then add to the $D(G)_{6s+1}$ sentences giving the connections including the triangle and pentagon
configuration.  We should also mark this connection as being made, so let $C_{3s} = C_{3s-1} \cup \{ (n_s, m_s) \}$. Otherwise, set $D(G)_{6s+1} = D(G)_{6s}$ and $C_{3s} =
C_{3s-1}$.

Step $6s+2$: In the next two steps we will make sure that for any atomic sentence $\phi$, either $\phi$ or $\lnot \phi$ will occur in $D(G)$. We find the first two elements $x$ and $y$ occurring in any sentence of $D(G)_{6s+1}$ for which neither $xGy$ nor $\lnot xGy$ is in the diagram so far. Now set
\[ \psi = \left\{ \begin{array}{ll} xGy & \mbox{if $yGx \in D(G)_{6s+1}$ or if $\{x, y\} = \{1, \nu_s(r)\}$ for some $r$}\\ \lnot xGy & \mbox{otherwise}\\
\end{array}
\right. \] and $D(G)_{6s+2} = D(G)_{6s+1} \cup \{\psi\}$.

Step $6s+3$: If $n_s = \nu_s(c)$ and $m_s = \nu_s(d)$ and $c = d \in D(T)$, then we'll have to say that $n_s = m_s$ in $G$, as well. In that case, set $D(G)_{6s+3} = D(G)_{6s+2} \cup \{ n_s = m_s \}$, and set $P_s = P_{s-1} \cup \{(n_s, m_s)\}$. Otherwise, set $D(G)_{6s+3} = D(G)_{6s+2} \cup \{ \lnot n_s =~ m_s \}$.

Step $6s+4$: Now we need two steps to be careful about the equalities we have just declared, to make sure that our construction is well-defined on isomorphism types. In particular, we will have to make sure that any pair that is equal has all the same connections. At this step, for all $(q,r) \in P_s$ and for all $t<s$, if $q T t$ and $(r, t) \notin C_{3s}$, then find the first available elements, and list the sentences necessary to form the triangle-pentagon pattern in them, as in step $6s+1$. Collect all such sentences as $J_s$, and set $D(G)_{6s+4} = D(G)_{6s+3} \cup J_s$ and $C_{3s+1} = C_{3s} \cup \{ (r,t) \}$.

Step $6s+5$: (This step is symmetric to $6s+4$.) For all $(q,r) \in P_s$ and for all $t<s$, if $r T t$ and $(q, t) \notin C_{3s}$, then find the first available numbers, and list the sentences necessary to form the triangle-pentagon pattern in them, as in step $6s+1$. Collect all such sentences as $J_s$, and set $D(G)_{6s+5} = D(G)_{6s+4} \cup J_s$ and $C_{3s+2} = C_{3s+1} \cup \{ (q,t) \}$.

Let $D(G) = \bigcup\limits_{s \in \omega} D(G)_s$. Now $D(G)$ is the diagram of the graph we wanted to build. If $T_1$ and $T_2$ are isomorphic directed graphs and we use this procedure to code them in undirected graphs $G_1$ and $G_2$ respectively, then an isomorphism $f: T_1 \rightarrow T_2$ induces an isomorphism from $G_1$ to $G_2$, which maps $\nu(x) \mapsto \nu(f(x))$, elements connecting $x$ and $y$ to those connecting $f(x)$ and $f(y)$, and which is otherwise the identity. Conversely, suppose that $G_1$ and $G_2$ are products of this procedure, perhaps the outputs given directed graphs $T_1$ and $T_2$, respectively, and that $G_1 \simeq G_2$. Then the isomorphism (we could call it $h$) must preserve the Archimedean point, and as a result must map vertex points to vertex points. Also, for any two vertex points $u_1$ and $u_2$, the subgraph $u_1$-triangle-pentagon-$u_2$ occurs in $G_1$ if and only if the sequence $h(u_1)$-triangle-pentagon-$h(u_2)$ occurs in $G_2$. Thus, there is also an isomorphism of the directed graphs $T_1$ and $T_2$. \end{proof}

\section{Arbitrary Computable Fields}

Intuition and experience tell us that the class of computable fields is quite complicated, perhaps so much that no classification could ever capture it. Previous work by Kudinov focused on existence of a computable ``Friedberg enumeration".

\begin{dfn} A Friedberg enumeration of $K$ up to isomorphism is a list of numbers, each of which is an index for a 
member of $K$, such that each isomorphism type from $K$ occurs exactly once in the list. The enumeration is said to be 
computable (or hyperarithmetical), when this list is.\end{dfn}

Goncharov and Knight had asked whether there was a computable Friedberg enumeration up to isomorphism of computable fields of fixed characteristic. Kudinov announced the following result:

\begin{thm}[Kudinov] There is no computable Friedberg enumeration of the computable fields of characteristic 0.\end{thm}

Knowing the complexity of the isomorphism problem for a class can tell us about the existence of Friedberg enumerations. The following is in \cite{gk}:

\begin{prop}[Goncharov-Knight] If $I(K)$ is hyperarithmetical and there is a hyperarithmetical Friedberg enumeration of $K$ up to isomorphism, then $E(K)$ is hyperarithmetical.\end{prop}

The idea of the proof is that if $E(K)$ is $\sum^1_1$ (which it must always be) and there is a hyperarithmetical Friedberg enumeration of $K$ up to isomorphism, then $E(K)$ is also $\prod^1_1$.

We can prove the following:

\begin{thm}\label{fields} When $K$ is the class of computable fields of some fixed characteristic, $E(K)$ is $\sum^1_1$ complete.\end{thm}

Then we have the following strengthening of Kudinov's result:

\begin{cor} For any $p$, prime or zero, there is no hyperarithmetical Friedberg enumeration up to isomorphism of computable fields of characteristic $p$.\end{cor}

\subsection{Borel Completeness for Fields: The Friedman-Stanley Embedding}

In 1989, Friedman and Stanley \cite{fs} showed that the class of countable fields of characteristic 0 is Borel complete, the maximal level of complexity in their sense. They proved this by constructing an embedding from graphs into fields (a Borel embedding is a Borel measurable function which is well-defined and injective on isomorphism types).

Friedman and Stanley assume that they are given a graph whose connectedness relation is $R$. From this they construct a field. We will use $\overline{F}$ to indicate the
algebraic closure of $F$, and $(S)$ for the smallest field containing $S$.

Consider $\{X_i\}_{i\in\omega}$, algebraically independent over $\mathbb{Q}$. Let $F_0$ be the composite of all of the $\overline{(\mathbb{Q}(X_i))}$, and define the
extension \[L(R) = F_0 (\{\sqrt{X_i + X_j} | iRj\})\] To deal with positive characteristic, we could replace $\mathbb{Q}$ with $\Fp$ and the square root with the $q$th root
where $q$ is relatively prime to $p$ both in this construction and throughout the following argument.  This function $L: Graphs \rightarrow Fields$ is both a Borel
measurable function under the usual product topology \cite{hjbook}, and well-defined on isomorphism classes.  The difficulty is in showing that

\begin{prop}\label{inj} L is injective on isomorphism classes.\end{prop}

\begin{proof}
In particular, it is difficult to show that $\sqrt{X_m + X_n}$ cannot be expressed as a rational function of the various $X_i$ and $\sqrt{X_j + X_k}$ where $\{j,k\} \neq
\{m,n\}$. The main difficulties appear when we consider only the composite of $\overline{(\mathbb{Q}(X_i))}$ and $\overline{(\mathbb{Q}(X_j))}$ (where $i \neq j$), and ask
whether it contains $\sqrt{X_i + X_j}$. I am grateful to W. Dwyer for the proof of the following lemma.  A different proof is given in the paper by Friedman and Stanley
\cite{fs}, and others still by Abhyankar \cite{abh} and Shapiro \cite{shapiro}.  In any case, proving the following lemma in positive characteristic is quite difficult.

\begin{lem} $\sqrt{X_i + X_j} \notin (\overline{\mathbb{Q}(X_i)} \cup \overline{\mathbb{Q}(X_j)})$, where $i \neq j$.\end{lem}

\begin{proof}A polynomial $p$ of degree $n$ in $\overline{(\mathbb{Q}(X_i))}$ gives a branched $n$-sheeted covering where the fiber over any point $a$ in
$\overline{(\mathbb{Q}(X_i))}$ is the set of roots of $p = a$, and branch points represent the multiple roots (a Riemann surface). A continuous function to find the roots
of the polynomial may be defined on this covering \textit{with branch points deleted}, but not in any neighborhood including the branch points themselves. We first consider
the possibility $$\sqrt{X_i + X_j} = \sum_{k = 1}^n a_kb_k$$ where $a_k$ is algebraic over $\mathbb{Q}(X_i)$, and $b_k$ is algebraic over $\mathbb{Q}(X_j)$. That is,
$\sum\limits_{k=1}^n a_k b_k$ gives one of the square roots.

Now for simplicity we can say that there is a single polynomial $p_i(Z)$ over $\mathbb{Q}(X_i)$ of which all the $a_k$ are roots, and similarly one $p_j$ over
$\mathbb{Q}(X_j)$ of which all the $b_k$ are roots. We can view the composite field $(\overline{\mathbb{Q}(X_i)}\cup \overline{\mathbb{Q}(X_j)})$ as
$\overline{\mathbb{Q}(X_i)} \times \overline{\mathbb{Q}(X_j)}$. Since $p_i$ and $p_j$ will each have only finitely many multiple roots and at most finitely many points at
which the coefficients are not defined, we can define continuous functions giving $a_k$ and $b_k$ on the relevant covering spaces of $$\overline{\mathbb{Q}(X_i)} \setminus
\{\mbox{these finitely many ``bad points" of }p_i\mbox{, say }e_t \}$$ and $$\overline{\mathbb{Q}(X_t)} \setminus \{\mbox{finitely many ``bad points" of }p_j\mbox{, say
}f_t \}$$ Thus the expression $\sum\limits_{k = 1}^n a_kb_k$ can be continuously defined on the relevant covering space of $\overline{\mathbb{Q}(X_i)} \setminus \{e_t\}
\times \overline{\mathbb{Q}(X_j)}\setminus \{f_t\})$. We can view this as a plane with finitely many vertical and horizontal lines deleted. Since the multiple roots of $Z^2
= X_i + X_j$ lie along the antidiagonal $X_i + X_j = 0$, there is clearly a neighborhood in which we can define $\sum\limits_{k = 1}^n a_kb_k$ as a continuous function, but
which contains points of the antidiagonal, so we cannot define $\sqrt{X_i + X_j}$ as a continuous function. Thus the two cannot be equal. To make the difference more
transparent, we could say that anywhere on this neighborhood we stay on the same branch of the right-hand side, but move from one branch to another on the left-hand side.

In the more general case that $$\sqrt{X_i + X_j} = \frac{\sum\limits_{k = 1}^n a_kb_k}{\sum\limits_{k = 1}^n c_kd_k}$$ we could simply write $$(\sqrt{X_i + X_j})(\sum_{k =
1}^n c_kd_k) = \sum_{k = 1}^n a_kb_k$$ and again the right-hand side can be defined continuously where the left-hand side cannot. \end{proof}

There does not seem to be a way to modify this proof to cover the positive characteristic case.  There is no apparent topology to replace the metric topology on the affine 
space, which we used here in declaring functions continuous or not.  Also, while we could talk about the number of values of the ``root function" and the right-hand-side 
``function," there are points in the proof at which it is not obvious that values will not collapse.

If we simply add more $X_i$, then more dimensions are added to the picture, but nothing really changes, since the diagonal for $\sqrt{X_i+X_j}$ is still in the same plane, and we can still find a neighborhood containing some point of it which contains no point of any line parallel to an axis. The next real problem comes up when we allow some square roots to be added. To simplify the task of visualization, and also to simplify the notation necessary, we will restrict the geometrical argument to a space whose $F_0$-dimension is the least possible to account for all $X_i$ used in the expression. This allows us to refer to codimension, allowing an economical way to describe the higher-dimensional generalizations of the fact that lines intersect in points, planes intersect in lines, and so forth.

\begin{lem} Let $F_0$ be as above. Then $$\sqrt{X_i + X_j} \notin F_0(\sqrt{X_j + X_k})$$ where $i$, $j$, and $k$ are distinct. \end{lem}

\begin{proof}Suppose not. First we will suppose again the simpler case where \[\sqrt{X_i+X_j} = \sum\limits_{s=1}^n \prod\limits_{t=1}^{m_s} a_{st}\] where each $a_{st}$ is algebraic over a single $X_q$ or is $f_{st} \sqrt{X_j + X_k}$ for some $f_{st} \in F_0$. Since there are only a finite number of such $X_q$ involved in the expression, let us collect, as before, polynomials $p_q$, one to account for all $a_{st}$ algebraic over a single $X_q$. The multiple roots of $p_q$ may be collected, as before, as $\{e_{q\gamma}\}$. Those used in the composition of $f_{st}$ may be collected as $\{d_{st\gamma}\}$. The multiple roots corresponding to $\sqrt{X_i + X_j}$ still form the diagonal (now a hyperplane, i.e. an algebraic surface of codimension 1) $X_i + X_j = 0$. Let $\hat{x}$ be a point of $X_i + X_j = 0$, and let $N$ be a ball around it of positive radius. Use $M_{\delta}$ to denote the (finitely many) hyperplanes $X_q = e_{q\gamma}$, $X_t = d_{st\gamma}$, and $X_j + X_k = 0$. Now $M^c_{\delta}$ (the complement of $M_{\delta}$) is open, so $N \cap (\bigcap\limits_{\delta} M_{\delta})$ is a neighborhood containing a point of $X_i + X_j = 0$ and no point of any $M_{\delta}$. Thus, there is a neighborhood in which we stay on a single branch of the right-hand side of the supposed equation, but cross a branch point of the left hand side.

Just as before, the extension to the more general case, $$\sqrt{X_i+X_j} = \frac{\sum\limits_{s=1}^n \prod\limits_{t=1}^{m_s} a_{st}}{\sum\limits_{s=1}^n \prod\limits_{t=1}^{m_s} b_{st}}$$ is quite easy. We clear the denominator and still have regions which are entirely fine for the right side of the equation but that the left finds unmanageable. \end{proof}

Alternately, we could consider a homomorphism $$F_0(\sqrt{X_j + X_k}) \rightarrow F_0(\sqrt{X_j + X_k})$$ which sends $X_k \mapsto 0$ but which is the identity on $\mathbb{Q} \cup \{X_{\ell}\}_{\ell \neq k}$. If the lemma failed, this homomorphism would show that $\sqrt{X_i + X_j} \in F_0$, contrary to the previous lemma. A similar alternate proof is possible for the following lemma.

Similarly, one can establish

\begin{lem} Let $F_0$ be as above. Then $\sqrt{X_i + X_j} \notin F_0(\{\sqrt{X_m + X_{n}}|\{m,n\} \neq \{i,j\}\})$.\end{lem}

\begin{proof}Suppose that the lemma fails. Then we suppose \[\sqrt{X_i+X_j} = \sum\limits_{s=1}^n \prod\limits_{t=1}^{m_s} a_{st}\] where each $a_{st}$ is algebraic over a single $X_q$ or is $f_{st} \sqrt{X_n+X_m}$ for some $f_{st} \in F_0$ and some $\{m,n\} \neq \{i,j\}$. Acting just as before, we denote by $p_q$ the polynomial accounting for all $a_{st}$ which are roots of some polynomial over $X_q$. The left-hand side of the equation still gives us multiple roots along the hyperplane $X_i + X_j = 0$, the roots of the $p_q$ still form hyperplanes parallel to the axes, just as before. The only difference from the previous case is that there are more hyperplanes of diagonal type ($X_m + X_n = 0$), but there is still a neighborhood in which the right side of the equation works, and the left does not. By this point, the usual extension to the more general form of a member of $$F_0(\{\sqrt{X_m + X_{n}}|\{m,n\} \neq \{i,j\}\})$$ is obvious. \end{proof}

We will need an additional fact. If two fields of this kind are isomorphic, the isomorphism will move $X_i$ to something interalgebraic with some $X_j$, since they are the elements whose algebraic closure is included. The change of $X_i$ to $X_j$ is clearly tolerable, since it merely corresponds to a permutation of the names of the vertices of a graph. However, we need to verify that the isomorphism does not foul up information on the connectedness relation.

\begin{lem} Let $A \sim B$ if and only if $\overline{\mathbb{Q}(A)} = \overline{\mathbb{Q}(B)}$. If $X \sim X_i$ and $Y \sim X_j$ then $\sqrt{X+Y} \notin (\overline{\mathbb{Q}(X_i)} \cup \overline{\mathbb{Q}(X_j)})$. \end{lem}

\begin{proof}The proof of this lemma is trivial. If $\sqrt{X+Y} \in (\overline{\mathbb{Q}(X_i)} \cup \overline{\mathbb{Q}(X_j)})$, then it is also in the (exactly equal) set $(\overline{\mathbb{Q}(X)} \cup \overline{\mathbb{Q}(Y)})$, in contradiction to the previous lemma. \end{proof}

\begin{lem} If $Y_i \sim X_i$ for all $i$, then $$\sqrt{Y_i + Y_j} \notin F_0(\{\sqrt{X_m + X_{n}}|\{m,n\} \neq \{i,j\}\})$$\end{lem}

\begin{proof}This proof is an almost equally obvious extension of earlier results. Let $F_1$ denote the composite field of all the $\overline{\mathbb{Q}(Y_i)}$. The field
$$F_0(\{\sqrt{X_m + X_{n}}|\{m,n\} \neq \{i,j\}\})$$ is equal to the field $$F_1(\{\sqrt{X_m + X_{n}}|\{m,n\} \neq \{i,j\}\})$$ We should note that each $X_i$ is the root
of some polynomial over $Y_i$. Now suppose \[\sqrt{Y_i + Y_j} = \sum\limits_{s=1}^n \prod\limits_{t=1}^{m_s} a_{st}\] where each $a_{st}$ is algebraic over a single $Y_i$
or is $f_{st} \sqrt{X_m + X_n}$ for some $f_{st} \in F_1$ and some $\{m,n\} \neq \{i,j\}$. Let $p_q$ again denote the polynomial accounting for all $a_{st}$ algebraic over
$Y_q$. The left-hand side still gives the same diagonal hyperplane $Y_i + Y_j = 0$. On the right-hand side, we have a finite union of hyperplanes parallel to the axes (the
multiple roots of $p_q$), and also some more exotic hypersurfaces. These hypersurfaces are those of the form $X_m + X_n = 0$. However, these are not equal to $Y_i + Y_j =
0$, so for each such hypersurface $P$ (of only finitely many) there is some neighborhood containing a point of $Y_i + Y_j = 0$ but no point of $P$.  Thus, we can still find
the necessary neighborhood in which the right-hand side of the equation is continuous and the left-hand side is not.  The more general element works as always. \end{proof}

We can now prove Proposition \ref{inj}. Suppose that $R$ and $S$ are two graphs, and that $L(R) \simeq L(S)$. Now by this isomorphism, each $X_i \in L(R)$ is mapped to some
$Y_i \sim X_j \in L(S)$. Certainly if $nRm$ then $\sqrt{X_n + X_m} \in L(R)$ and thus $\sqrt{Y_n + Y_m} \in L(S)$. By the last lemma, if $Y_n \sim X_p$ and $Y_M \sim X_q$,
the last statement implies that $\sqrt{X_p + X_q} \in L(S)$, so by the previous lemma $pSq$ (since $\sqrt{X_p + X_q} \notin F_0(\{\sqrt{X_m + X_{n}}|\{m,n\} \neq
\{i,j\}\})$ (that is, $X_p + X_q$ only had a square root if we put one in to account for a connection of $p$ and $q$ in $S$). Similarly, we can argue that if $nSm$, then
the corresponding elements are connected in $R$. Thus $R \simeq S$. \end{proof}

\subsection{Computable Construction of the Friedman-Stanley Embedding}

It will turn out that a similar embedding produces computable fields from computable graphs, amounting to a reduction $E($Graphs$) \leq_1 E($Fields$)$. This will complete
the proof. The only real modification necessary is to guarantee that if we start with a computable graph, we end up with a computable field. Important background work on
computable fields may be found in \cite{frsh} \cite{mnfields} \cite{mr}.

We should note that since the class of fields of given characteristic has $\prod^0_2$ axioms (stating that it is a commutative ring, plus the condition that for any element
there exists a multiplicative inverse), $I(K)$ is $\prod^0_2$. Given a computable directed graph with connectedness relation $R$, consider $\{X_i\}_{i\in\omega}$,
algebraically independent over $\mathbb{Q}$. Let $G_0$ be a computable field isomorphic to the composite of all of the $\overline{\mathbb{Q}(X_i)}$, and let $M(R)$ be the
extension $G_0 (\{\sqrt{X_i + X_j} | iRj\})$. It remains to verify

\begin{prop}$M(R)$ has a computable copy, whose index can be obtained effectively from that of $R$.\end{prop}

\begin{proof}Consider the language of fields, plus countably many constants $X_i$, with the theory of algebraically closed fields of characteristic 0 and the sentences
stating that the $X_i$ are algebraically independent. This theory is complete and decidable (since the theory of algebraically closed fields alone proves quantifier
elimination), and so it has a computable model. Call this computable model $G$.

Once we have $G$, there is a c.e.\ $G_0^* \subseteq G$ which contains exactly those members of $G$ which are algebraic over a single $X_i$. Further, there is a c.e.\ $R^*
\subseteq G$ consisting of exactly $\{\sqrt{X_i + X_j} | iRj\}$. With these two sets, we can enumerate the elements of the smallest subfield containing $G_0^* \cup R^*$,
and we will call this $F^* \subseteq G$. Note that $F^*$ has c.e.\ universe. Let $e(R)$ be the index of the function with which we enumerate $D(F^*)$, and note that we can
find it effectively in a uniform way from an index for $D(R)$. Now by padding, we can replace the c.e.\ field $F^*$ by a field whose universe is computable. Let $\tilde{M}$
be the set of all $(a,s)$ where $a \in F^*_s \setminus F^*_{s-1}$. It is clear that the reduct of $\tilde{M}$, with the operations $(a,s) +_M (b,t) := (a +_G b, r)$ and
$(a,s) \cdot_M (b,t) := (a \cdot_G b, r)$ is a computable field isomorphic to $M(R)$, where in each case $r$ is the least such that the desired elements are in $F^*_r$. It
is also clear that an index for the function with which we enumerate this field is effectively obtained from $e(R)$ in a uniform way. This completes the proof both of the
proposition and of Theorem \ref{fields} in the case of characteristic zero. \end{proof}

\section{Algebraically Closed Fields and Vector Spaces}

While arbitrary computable fields could be expected to be difficult to classify, the algebraically closed fields ought to be much simpler, since isomorphism can be checked
by just comparing transcendence degree. Vector spaces, interesting in their own right, give us a glimpse of the methods to be used on algebraically closed fields while
using simpler notation and being more intuitively accessible.

Vector spaces are classified by their dimension. A bijection of bases of vector spaces induces an isomorphism of the spaces. Since it is computationally not too hard to find a basis \cite{mnvs}, we should expect that the isomorphism problem will be rather simple.

\begin{thm} Fix an infinite computable field $F$, and let $K$ be the class of computable vector spaces over $F$. 
Then $E(K)$ is $\Pi^0_3$ complete. \end{thm}

\begin{proof}To see that the problem is $\Pi^0_3$, we first note that $I(K)$ is $\Pi^0_2$, since the class can be
axiomatized with $\Pi^0_1$ axioms declaring that the structure is an additive Abelian group and a module, and since
the statement ``$a$ is the index of some structure" is $\Pi^0_2$. We will formulate a sentence to say that
$\mathcal{A}_a$ and $\mathcal{A}_b$ have the same dimension.  The statement \[D_n = \exists (x_1, ..., x_n)  
\bigwedge_{\bar{a}\neq 0}\hspace{-.15in}\bigwedge (\sum_{i=1}^n a_ix_i \neq 0)\] (which is $\Sigma^0_2$) states that
there are at least $n$ linearly independent elements. Now we can write \[\forall n [\mathcal{A}_a \models D_n
\leftrightarrow \mathcal{A}_b \models D_n]\] (which is $\Pi^0_3$), and the conjunction of this statement with $a,b
\in I(K)$ defines $E(a,b)$.

To see that the problem is $\Pi^0_3$ complete, we will take a known $\Pi^0_3$ complete set $S$, and construct a 
uniformly computable sequence of vector spaces $(V(n))_{n\in\omega}$ such that $V(n) \simeq V^\infty$ if and only if 
$n\in S$, where $V^\infty$ is an infinite dimensional vector space. This will show that if $e$ is the index of 
$V^{\infty}$, $E(a,e)$ is $\Pi^0_3$ complete.

Define $Cof = \{e | W_e$ is cofinite$\}$, where $W_e$ is the domain of the $e$th partial recursive function, and take 
$S = \overline{Cof}$. Now $V(n)$ will be a vector space with dimension the same as $card (\overline{W_n})$. Note that 
$Cof$ is $\Pi^0_3$ complete \cite{soare}.

We begin by letting $\{a_i\}_{i\in\omega}$ be a basis for a computable copy of $V^{\infty}$, and by considering a set $\{b_i\}_{i\in\omega}$. Let $\lambda_i$ enumerate the
linear combinations of $b_i$. We will construct $f$ and $D(V(n))$ to meet the following conditions:

\smallskip
\begin{tabular}{rl} $P^1_e$ : & $\lambda_e = 0 \in D(V(n))$ or $\lnot \lambda_e = 0 \in D(V(n))$\\
$P^2_e$ : & $b_e \in dom(f)$\\ 
$Q_e$ : & $f(b_e)$ is linearly independent of $b_0, ... b_{e-1}$ if and only \\
        & if $e \notin W_n$\\ 
\end{tabular}
\smallskip

Thus, $im(f)$ will be a subspace of $V^\infty$ which we will call $V(n)$, and will have the desired properties. We say that $P^1_e$ requires attention at stage $s$ if $\lambda_e = 0 \notin D(V(n))_s$ and $\lnot \lambda_e = 0 \notin D(V(n))_s$ and $i < \frac{s}{2}$ wherever $b_i$ is included in $\lambda_e$. $Q_e$ requires attention at stage $s$ if $e \in W_{n,s}$ and $Q_e$ has not received attention previously. Let $f_0$ be the homomorphism induced by taking $b_i \mapsto a_i$, and $D(V(n))_0 = \emptyset$.

At stage $2s$, let $\ell_s$ be the first $e$ such that $P^1_e$ requires attention. We can then find whether $f_s(\lambda_{\ell_s}) = 0$. If so, let $D(V(n))_{2s+1} = D(V(n))_{2s} \cup \{\lambda_{\ell_s} = 0\}$. Otherwise, add $\lnot \lambda_{\ell_s} = 0$.

At stage $2s+1$, we will work on $Q_e$. Let $x_s$ be the least $e$ such that $Q_e$ requires attention. Let $\delta_s(b_1, ..., b_N)$ be the (finite) conjunction of all
sentences in $D(V(n))_{2s+1}$. Because the theory of vector spaces is strongly minimal and the closure of any set except zero is infinite, there is some $\hat{a}_{x_s} \in
Cl(a_0, a_1, ..., a_{x_s-1})$ such that $V^\infty \models \delta_s(f_s(b_0), f_s(b_1), ..., f_s(b_{x_s-1}),\newline \hat{a}_{x_s}, f_s(b_{x_s+1}), ..., f_s(b_N))$. We then
define \[ f_{s+1}(b_i) = \left\{ \begin{array}{ll} \hat{a}_{x_s} & \mbox{if $i=x_s$}\\ f_{s}(b_i) & \mbox{otherwise}\\ \end{array} \right. \] and let $D(V(n))_{2s+2} =
D(V(n))_{2s+1}$.

Now it is clear that for any $e$, there is some stage at which $P^1_e$ will receive attention, after which it will never again be injured. Further, $P^2_e$ is never injured. Finally, $Q_e$ will either be satisfied at every stage ($e \notin W_n$), or it will receive attention at some stage, after which it will never be injured again. \end{proof}

\begin{thm} The isomorphism problem for the class of computable algebraically closed fields of fixed 
characteristic is $\Pi^0_3$ complete. \end{thm}

\begin{proof} The construction is very similar to the previous one, except that dimension is replaced by transcendence degree, and $F$-linear combinations are replaced by
polynomial combinations.  The sentence \[T_n = \exists (x_0, ... x_n)  \bigwedge\hspace{-.1in}\bigwedge 
(\sum_{i=1}^{2^n} (\prod_{j=0}^{n} x_{j}^{p_{ij}}) \neq 0)\] where
the conjunction is taken over all $\{p_{ij}\}_{1 \leq i \leq 2^n, 1 \leq j \leq n}$ states that the field has transcendence degree at least $n$.  Note that even in 
positive characteristic, this conjunction is necessarily infinite, since there is at least one irreducible polynomial in every degree.  Now we can write \[\forall
n [\mathcal{A}_a \models T_n \leftrightarrow \mathcal{A}_b \models T_n]\] in conjunction with $a,b \in I(K)$ for 
$E(a,b)$. Note that $I(K)$ is $\Pi^0_2$.

We will build a uniformly computable sequence $(F(n))_{n \in \omega}$ of algebraically closed fields such that $F(n)$ will have transcendence degree equal to the
cardinality of $\overline{W_n}$. Just as before, this will establish the theorem.

We begin by letting $\{a_i\}_{i\in\omega}$ be a transcendence base for a computable copy of $F^{\infty}$ (the algebraically closed field of infinite transcendence degree), and by considering a set $\{b_i\}_{i\in\omega}$. Let $\lambda_i$ enumerate the polynomial combinations of $b_i$. We will construct $f$ and $D(F(n))$ to meet the following conditions:

\smallskip
\begin{tabular}{rl} $P^1_e$ : & $\lambda_e = 0 \in D(F(n))$ or $\lnot \lambda_e = 0 \in D(F(n))$\\
$P^2_e$ : & $b_e \in dom(f)$\\ $Q_e$ : & $f(b_e)$ is algebraically independent of $b_0, ... b_{e-1}$ if and only if we \\ 
          & have $e \notin W_n$\\ 
\end{tabular}
\smallskip

Thus, $im(f)$ will be a subset of $F^\infty$ which we will call $F(n)$, and will have the desired properties. We say that $P^1_e$ requires attention at stage $s$ if 
$\lambda_e = 0 \notin D(F(n))_s$ and $\lnot \lambda_e = 0 \notin D(F(n))_s$ and $i < \frac{s}{2}$ wherever $b_i$ is included in $\lambda_e$. $Q_e$ requires attention at stage $s$ if $e \in W_{n,s}$ and $Q_e$ has not received attention previously. Let $f_0$ be the homomorphism induced by taking $b_i \mapsto a_i$, and $D(F(n))_0 = \emptyset$.

At stage $2s$, let $\ell_s$ be the first $e$ such that $P^1_e$ requires attention. We can then find whether $f_s(\lambda_{\ell_s}) = 0$. If so, let $D(F(n))_{2s+1} = D(F(n))_{2s} \cup \{\lambda_{\ell_s} = 0\}$. Otherwise, add $\lnot \lambda_{\ell_s} = 0$.

At stage $2s+1$, we will work on $Q_e$. Let $x_s$ be the least $e$ such that $Q_e$ requires attention. Let $\delta_s(\overline{b})$ be the (finite) conjunction of all
sentences in $D(F(n))_{2s+1}$. Again using strong minimality, there is some $\hat{a}_{x_s}$ in the closure of the set 
of $a_i$ for $i<X_s$ such that \[F^\infty \models
\delta_s(f_s(b_0), f_s(b_1), ..., f_s(b_{x_s-1}), \hat{a}_{x_s}, f_s(b_{x_s+1}), ..., f_s(b_N))\]  We then define

\[ f_{s+1}(b_i) = \left\{ \begin{array}{ll} \hat{a}_{x_s} & \mbox{if $i=x_s$}\\ f_{s}(b_i) & \mbox{otherwise}\\
\end{array}
\right. \] and let $D(F(n))_{2s+2} = D(F(n))_{2s+1}$.

Again it is clear that for any $e$, there is some stage at which $P^1_e$ will receive attention, after which it will never again be injured. Further, $P^2_e$ is never injured. Finally, $Q_e$ will either be satisfied at every stage ($e \notin W_n$), or it will receive attention at some stage, after which it will never be injured again. \end{proof}

\section{Real Closed Fields}

Knowing that algebraically closed fields are very simple, but having the suspicion that fields in general will be
more complicated, it would be interesting to know about other restrictions of fields. K.\ Manders suggested the
example of real closed fields, whose model theory is reasonably well-behaved, but which is unstable. The complexity is
all in the infinite elements.

\begin{thm} If $K$ is the class of Archimedean real closed fields, then $E(K)$ is $\Pi^0_3$ complete.\end{thm}

\begin{proof}Recall that a real closed field is an ordered field satisfying the additional condition that each odd-degree polynomial has a root. Thus, the class of real
closed fields can be axiomatized by a computable infinitary $\Pi_2$ sentence, as can the class of Archimedean real 
closed fields (by adding the sentence that for each
element $x$, some finite multiple of $1$ is greater than $x$). Archimedean real closed fields are classified simply by the cuts that are filled, so the statement \[\forall
x,z\hspace{0.05in} \exists \hat{x}, \hat{z} \bigwedge_{q \in \mathbb{Q}}\hspace{-.15in}\bigwedge [(\mathcal{A}_a \models q \leq x \Leftrightarrow \mathcal{A}_b \models q
\leq \hat{x})]\wedge\] \nopagebreak \[\wedge \bigwedge_{q \in \mathbb{Q}}\hspace{-.15in}\bigwedge [(\mathcal{A}_b \models q \leq z \Leftrightarrow \mathcal{A}_a \models q
\leq \hat{z})]\] defines the relation $\mathcal{A}_a \simeq \mathcal{A}_b$, showing that it is, at worst, $\Pi^0_3$.

\begin{lem} \label{indepreal} There exists a uniformly computable sequence $\{a_i\}_{i \in \omega}$ of real numbers 
which are algebraically independent.\end{lem}

\begin{proof} Lindemann's theorem, a known result, states (in one form) that if $\lambda_1, \lambda_2, \dots,
\lambda_k$ are algebraic numbers linearly independent over the rationals, then $e^{\lambda_1}, \dots, e^{\lambda_k}$
are algebraically independent \cite{baker}.  Further, it is well known that the set $\{ \sqrt{2}, \sqrt{3}, \sqrt{5},
\dots \}$ is linearly independent (a proof may be found in 
\cite{besicovitch}). 
\end{proof}

Now consider the language $(+,\cdot, 0, 1, \{a_i\}_{i \in \omega})$, and the theory consisting of the axioms of real closed fields and the set $A$. This is a complete decidable theory, and so has a computable model. The set of elements algebraic over the set of $a_i$ is c.e.\ and by padding we can find a computable structure $\mathcal{M} = (M, +, \cdot, 0, 1, \leq, \{a_i\}_{i \in \omega})$ where $M$ is the real closure of $\{a_i\}_{i \in \omega}$.

We are now ready to prove the completeness part of the theorem. We will start with a set $\{b_{ij}\}_{i, j \in 
\omega}$, $D_{-1} = \emptyset$, a list $\lambda_i$ of all
(positive) atomic sentences in the language of ordered fields with constants $b_{ij}$, and a function $f_{-1}: b_{ij} \mapsto a_i$. Let $S = \{n | \hspace{0.05in} \forall i
\exists j \forall z R(i,j,z,n)\}$ be an arbitrary $\Pi^0_3$ set. The proof will be similar to the proof for 
algebraically closed fields in that we will construct a
uniformly computable sequence $F(n)$ of real closed fields such that $F(n) \simeq \mathcal{M}$ exactly when $n \in S$. Set $m_{t,s} = 0$ and $I_{s,-1} = D_{-1} = \emptyset$
for all $(t,s)$. We wish to meet the following requirements:

\smallskip
\begin{tabular}{rl} $P^1_e$ : & $\lambda_e = 0 \in D$ or $\lnot \lambda_e = 0 \in D$\\
$P^2_e$ : & $b_{i,j} \in dom(f)$\\ $Q_e$ : & For each $i$, if $n \in S$ there exists some $b_{ij}$ which fills the same\\ & cut in $F(n)$ that $a_i$ fills in $M$. Otherwise there is no such $b_{ij}$.\\ \end{tabular}
\smallskip

At stage $2s$, find the least $i$ such that $\lambda_i \notin D_{s-1}$ and $\lnot \lambda_i \notin D_{s-1}$. Without loss of generality say that $\lambda_i$ is of the form $p(b_{j_0}, \dots, b_{j_n}) Q 0$ where $Q$ is $=$ or $\leq$. Now if $\mathcal{M} \models p(f_{s-1}(b_{j_0}), \dots, f_{s-1}(b_{j_n})) Q 0$ then set $D_s = D_{s-1} \cup \{\lambda_i\}$. Otherwise, $D_s = D_{s-1} \cup \{\lnot \lambda_i\}$.

At stage $2s+1$, for each $t<s$ we will check whether for all $z<s$ we have $R(t, m_{t,s}, z, n)$. If this holds, set $I_{s,t} = I_{s,(t-1)}$ and $m_{t,(s+1)} 
= m_{t,s}$. Otherwise
set $m_{t,(s+1)} = \min \{k>m_{t,s} | b_{tk} \notin \bigcup_{q<s} I_q$ and $b_{tk} = b_{tm_{t,s}} \notin D_s\}$. Now let $\delta_s (\vec{b},x)$ be the conjunction of all
sentences in $D_s$ true on $b_{tm_{t,s}}$, and note that the set it defines contains an interval, since everything in it must be true in $\mathcal{M}$ of $a_t$. There is
some rational $r_{t,m_{t,s}}$ in this interval, which we can find effectively, and we set $I_{s,t} = I_{s,(t-1)} \cup \{(t, m_{t,s}, r_{t,m_{t,s}}\} \cup \{(t, k,
r_{t,m_{t,s}}) | D_s \vdash b_{tk} = b_{tm_{t,s}}\}$. Note that this last addition can be made effectively, since only finitely many $b_{tk}$ will have been mentioned, and
since the theory of real closed fields is decidable. Let $I_s = \bigcup_{t<s} I_{s,t}$. Now we change the function:

\[ f_{s}(b_{ij}) = \left\{ \begin{array}{ll} \hat{r} & \mbox{if $(i, j, \hat{r}) \in I_s$}\\ f_{s-1}(b_{ij}) & \mbox{otherwise}\\
\end{array}
\right. \]

Let $D(F(n)) = \bigcup_{s \in \omega} D_s$, and we will call the structure whose diagram this is $F(n)$. Now if $n \in S$, for each $i$ there will eventually be some $j$ such that we always leave $b_{ij} \mapsto a_i$, so the cut corresponding to $a_i$ is filled, and $F(n) \simeq \mathcal{M}$. Otherwise, for some $i$, each $b_{ij}$ will be mapped to some rational, and the isomorphism will fail. \end{proof}

When we add positive infinite elements, however, we have a great deal of freedom in the structure of the field.

\begin{thm} \label{rcf} If $K$ is the class of real closed fields, then $E(K)$ is $\Sigma^1_1$ 
complete.\end{thm}

\begin{proof}We will say that $a \preceq b$ exactly when $a^n \leq b$ for 
some $n \in \omega$. We say that \mbox{$a
\approx b$} ($a$ is comparable to $b$) if $a \preceq b$ and $b \preceq a$. The proof will depend on realizing an
arbitrary computable linear order as the order type of the comparability classes of infinite elements.

\begin{lem} Given a computable linear order $L$, there is a computable structure $R(L)^* = (\bar{R}, +, \cdot, 0, 1, \leq, \{X_i\}_{i \in \omega})$, an expansion of a real closed field, in which $X_i \preceq X_j$ if and only if $i \leq_L j$, whose index is computable from an index for $L$.\end{lem}

Consider the language of ordered fields, plus infinitely many constants $X_i$, with the theory of real closed fields, and the sentences for each $i$ stating that $X_i$ is greater than any polynomial in $\{X_j | j <_L i\}$, and that all are greater than polynomials in $1$. This is a complete, decidable theory, and thus has a computable model $G$. There is a c.e.\ subset $\tilde{R} \subseteq G$ containing exactly the elements algebraic over $\{X_i\}_{i \in \omega}$. From an index for $L$, we can effectively find an index $e(L)$ for the function enumerating $\tilde{R}$. Again we can pad to find an isomorphic structure $R(L)^*$ with computable universe, as claimed. Let $R(L)$ denote the reduct of $R(L)^*$ to the language of ordered fields.

So we have encoded arbitrary linear orders into real closed fields, and all that remains is to make sure that this
operation is well-defined and injective on isomorphism types. The well-definedness is clear, since an isomorphism of
linear orders would just amount to a permutation of the labels for the $X_j$. It is also clear that if $h$ is an
isomorphism $h: R(L_1) \rightarrow R(L_2)$, then for $a, b \in R(L_1)$, $a \preceq b$ if and only if $h(a) \preceq
h(b)$, but it requires some verification to see that for $a$ in the comparability class of some $X_i$, $h(a)$ must be
in the comparability class of some $X_j$. Once this is shown, $h$ will induce an isomorphism of orders $\tilde{h}$,
where if $h$ maps the class of $X_i$ to that of $X_j$, then $\tilde{h}: i \mapsto j$. I am grateful to L.\ van den
Dries for suggesting the proof of the following lemma.

\begin{lem} Let $F$ be $R(L)$ for some linear order $L$. Let $C$ be a positive infinite comparability class of
elements of $F$. Then $C$ is the comparability class of one of the $X_i$.\end{lem}

\begin{proof} Suppose we have a real closed field $K$, and we add a single positive infinite element $x > K$. Let
$K((X^{\mathbb{Q}}))$ denote the set of formal series $f = \sum\limits_{q \in \mathbb{Q}} a_q X^q$, where $a_q \in K$
and $a_q = 0$ except for $q$ in some well-ordered set. There is an isomorphism $rcl(K(x)) \simeq K((X^{\mathbb{Q}}))$
mapping $x \mapsto X^{-1}$. Now suppose that $y \in K((X^{\mathbb{Q}}))$, and $y = \sum\limits_{q \in \mathbb{Q}} b_q
X^q$. Further, suppose that for all $x \in K$, we have $y > x$ (that is, $y$ is an infinite element over $K$).  Let
$\hat{q}$ be the least such that $b_{\hat{q}} \neq 0$.  I claim that $\hat{q} < 0$.  If $\hat{q} \geq 0$, then $y \leq
b_0 +1$, but $b_0+1 \in K$, giving a contradiction.  So $\hat{q}<0$.  Now $X^{\hat{q}-1} >y$; that is,
$(X^{-1})^{1-\hat{q}} > y$, so $y \approx t^{-1}$.  Thus, $rcl(K(x))$ has exactly one more comparability class than
$K$.

Given this, the lemma is relatively easy. Using the previous paragraph as an induction step, it is easy to show that
for $L$ a finite linear order, the lemma holds. Further, since any element in $R(L)$ is algebraic over finitely many
$X_i$, $R(L) = \bigcup\limits_I rcl(R(\{X_i\}_{i \in I}))$ where $I$ is a finite subset of $\omega$. This completes
the proof of both the lemma and the theorem. \end{proof}\end{proof}

It is worthwhile to note that Theorem \ref{rcf} implies the characteristic $0$ case of Theorem \ref{fields}.  This 
proof is certainly simpler.  However, the earlier proof covers positive characteristic and stresses the relationship 
with Borel complexity.  Also, it offers an opportunity to simplify, at least for positive characteristic, the 
difficult argument of the Friedman--Stanley paper.

\bibliographystyle{amsplain}
\bibliography{ecl}

\providecommand{\bysame}{\leavevmode\hbox to3em{\hrulefill}\thinspace}
\begin{thebibliography}{10}

\bibitem{abh}
S.~Abhyankar, \emph{On the compositum of algebraically closed subfields},
  Proceedings of the American Mathematical Society \textbf{7} (1956), 905--907.

\bibitem{baker}
A.~Baker, \emph{Transcendental number theory}, Cambridge University Press,
  1975.

\bibitem{bkbook}
H.~Becker and A.~S. Kechris, \emph{The descriptive set theory of polish group
  actions}, London Mathematical Society Lecture Note Series, vol. 232,
  Cambridge University Press, 1996.

\bibitem{besicovitch}
A.~S. Besicovitch, \emph{On the linear independence of fractional powers of
  integers}, Journal of the London Mathematical Society \textbf{15} (1940),
  3--6.

\bibitem{ecl2}
W.~Calvert, \emph{The isomorphism problem for classes of computable {A}belian
  groups}, preprint, 2003.

\bibitem{fs}
H.~Friedman and L.~Stanley, \emph{A {B}orel reducibility theory for classes of
  countable structures}, Journal of Symbolic Logic \textbf{54} (1989),
  894--914.

\bibitem{frsh}
A.~Frolich and J.~C. Sheperdson, \emph{Effective procedures in field theory},
  Philosophical Transactions of the Royal Society of London, Series A
  \textbf{248} (1956), 407--432.

\bibitem{gk}
S.~S. Goncharov and J.~F. Knight, \emph{Computable structure and non-structure
  theorems}, Algebra and Logic \textbf{41} (2002), 351--373.

\bibitem{hkss}
D.~Hirschfeldt, B.~Khoussainov, R.~Shore, and A.~M. Slinko, \emph{Degree
  spectra and computable dimensions in algebraic structures}, Annals of Pure
  and Applied Logic \textbf{115} (2002), 71--113.

\bibitem{hjbook}
G.~Hjorth, \emph{Classification and orbit equivalence relations}, American
  Mathematical Society, 1999.

\bibitem{hjkelo}
G.~Hjorth, A.~S. Kechris, and A.~Louveau, \emph{Borel equivalence relations
  induced by actions of the symmetric group}, Annals of Pure and Applied Logic
  \textbf{92} (1998), 63--112.

\bibitem{hodges}
W.~Hodges, \emph{What is a structure theory?}, Bulletin of the London
  Mathematical Society \textbf{19} (1987), 209--237.

\bibitem{mnvs}
G.~Metakides and A.~Nerode, \emph{Recursively enumerable vector spaces}, Annals
  of Mathematical Logic \textbf{11} (1977), 146--171.

\bibitem{mnfields}
\bysame, \emph{Effective content of field theory}, Annals of Mathematical Logic
  \textbf{17} (1979), 289--320.

\bibitem{moro}
A.~S. Morozov, \emph{Functional trees and automorphisms of models}, Algebra and
  Logic \textbf{32} (1993), 28--38.

\bibitem{nies}
A.~Nies, \emph{Undecidable fragments of elementary theories}, Algebra
  Universalis \textbf{35} (1996), 8--33.

\bibitem{mr}
M.~O. Rabin, \emph{Computable algebra, general theory and theory of computable
  fields}, Transactions of the American Mathematical Society \textbf{95}
  (1960), 341--360.

\bibitem{rabsco}
M.~O. Rabin and D.~Scott, \emph{The undecidability of some simple theories},
  preprint.

\bibitem{shapiro}
D.~B. Shapiro, \emph{Composites of algebraically closed fields}, Journal of
  Algebra \textbf{130} (1990), 176--190.

\bibitem{shelah}
S.~Shelah, \emph{Classification of first order theories which have a structure
  theorem}, Bulletin of the American Mathematical Society (New Series)
  \textbf{12} (1985), 227--232.

\bibitem{soare}
R.~I. Soare, \emph{Recursively enumerable sets and degrees}, Springer-Verlag,
  1987.

\end{thebibliography}

\end{document}